\documentclass[12pt]{amsart}
\usepackage{verbatim}
\usepackage{amssymb,amsmath}

\newtheorem{lemma}{Lemma} [section]

\newtheorem{cor}[lemma]{Corollary}
\newtheorem{prop}[lemma]{Proposition}

\theoremstyle{remark}
\newtheorem*{remark}{Remark}

\newcommand{\F}{{\mathbb F}}

\newcommand{\Z}{{\mathbb Z}}

\newcommand{\CC}{{\mathcal C}}

\newcommand{\tth}{^{\operatorname{th}}}


\begin{document}



\title[Permutation polynomials]{Some families of permutation polynomials
over finite fields}

\author{Michael E. Zieve}
\address{Center for Communications Research, 805 Bunn Drive, Princeton NJ 08540}

\email{zieve@math.rutgers.edu}
\urladdr{www.math.rutgers.edu/$\sim$zieve/}

\begin{abstract}
We give necessary and sufficient conditions for
a polynomial of the form $x^r(1+x^v+x^{2v}+\dots+x^{kv})^t$
to permute the elements of the finite field $\F_q$.
Our results yield especially simple criteria in case 
$(q-1)/\gcd(q-1,v)$ is a small prime.
\end{abstract}


\maketitle


\section{Introduction}

A polynomial over a finite field is called a \emph{permutation polynomial}
if it permutes the elements of the field.  These polynomials first arose
in work of Betti~\cite{B}, Mathieu~\cite{M} and Hermite~\cite{H} as a way
to represent permutations.
A general theory was developed by Hermite~\cite{H} and Dickson~\cite{D}, with
many subsequent developments by Carlitz and others.

It is a challenging problem to produce permutation polynomials of `nice' forms.
Recently, Akbary, Wang and Wang \cite{AW0,Wang} studied
binomials of the form $x^u+x^r$ over $\F_q$ in the case
that $d:=\gcd(q-1,u-r)$ satisfies $(q-1)/d\in\{3,5,7\}$.  Their results were
surprising: they gave necessary and sufficient criteria for such binomials to
permute $\F_q$, in terms of the period of a (generalized) Lucas sequence
in $\F_q$.  Their proofs were quite complicated, using lengthy calculations
involving coefficients of Chebychev polynomials, lacunary sums of binomial
coefficients, determinants of circulant matrices, and various unpublished
results about factorizations of Chebychev polynomials, among other things.
Also, their proofs required completely different arguments in each of the
cases $(q-1)/d\in\{3,5,7\}$.

One naturally wonders whether there might be a uniform approach which
works for arbitrary $d$, and yields the results of \cite{Wang,AW0} as
special cases.  We
present such an approach in this paper, giving short and simple proofs
which do not use any of the above-mentioned ingredients.
Our results apply to the more general class
of polynomials $f(x):=x^r h_k(x^v)^t$,
where $h_k(x):=x^{k-1}+x^{k-2}+\dots+1$ and $r,v,k,t$ are positive integers.
The forthcoming paper \cite{AAW} uses the same methods as \cite{AW0} to
prove some partial results in case $t=1$ and $v\mid (q-1)$.

The statements of our results use the notation
$s:=\gcd(v,q-1)$, $d:=(q-1)/s$, and $e:=v/s$.  Note that $\gcd(d,e)=1$.
Also $\mu_d$ denotes the set of $d\tth$ roots of unity in $\F_q$,
and $p$ is the characteristic of~$\F_q$.

Our first result gives necessary and sufficient conditions for $f$ to
be a permutation polynomial:

\begin{prop}
\label{intropl}
$f$ permutes\/ $\F_q$ if and only if all of the following hold:
\begin{enumerate}
\item $\gcd(r,s)=\gcd(d,k)=1$
\item $\gcd(d,2r+vt(k-1))\le 2$
\item $k^{st}\equiv (-1)^{(d+1)(r+1)}\pmod{p}$
\item $g(x):=x^r((1-x^{ek})/(1-x^e))^{st}$ is injective on
$\mu_d\setminus\mu_1$
\item $(-1)^{(d+1)(r+1)} \notin g(\mu_d\setminus\mu_1)$.
\end{enumerate}
\end{prop}

In case $d$ is an odd prime, this specializes to the following:

\begin{cor}
\label{introneat}
Suppose the first three conditions of Proposition~\ref{intropl} hold,
and $d$ is an odd prime.  Pick $\omega\in\F_q$ of order~$d$.
Then $f$ permutes $\F_q$ if and only if there exists $\theta\in\F_d[x]$
with $\theta(0)=0$ such that $(2r+(k-1)vt)x+\theta(x^2)$ permutes $\F_d$
and, for every $i$ with $0<i<d/2$, we have
\[
\omega^{\theta(i^2)} = \left(\frac{\omega^{ike}-\omega^{-ike}}{\omega^{ie}-
\omega^{-ie}}\right)^{st}.
\]
\end{cor}

In the cases $d=3,5,7$ studied in \cite{Wang} and \cite{AW0}, it remains
to consider permutation polynomials of $\F_d$ of certain forms.  This is
quite simple to analyze directly, and it is also a consequence of the results
of Betti (1851) and Hermite (1863).  The conclusion is as follows:

\begin{cor}
Suppose the first three conditions of Proposition~\ref{intropl} hold,
and $d$ is an odd prime.  Pick $\omega\in\F_q$ of order~$d$.
\begin{enumerate}
\item If
\begin{equation}
\tag{$*$}
\frac{\zeta^k-\zeta^{-k}}{\zeta-\zeta^{-1}}\in\mu_{st} \text{ for every
$\zeta\in\mu_d\setminus\mu_1$}
\end{equation}
then $f$ permutes\/ $\F_q$.
\item If $d=3$ then $f$ always permutes\/ $\F_q$.
\item If $d=5$ then $f$ permutes\/ $\F_q$ if and only if $(*)$ holds.
\item If $d=7$ then $f$ permutes\/ $\F_q$ if and only if either $(*)$ holds or
there exists $\epsilon\in\{1,-1\}$ such that
\[
\left(\frac{\omega^{ike}-\omega^{-ike}}{\omega^{ie}-
\omega^{-ie}}\right)^{st} = \omega^{2\epsilon (2r+(k-1)vt)i}
\]
for every $i\in\{1,2,4\}$.
\end{enumerate}
\end{cor}

It is straightforward to deduce the results of \cite{Wang,AW0,AAW}
from this result, by writing the generalized Lucas sequences
in terms of roots of unity.  However, our formulation seems to be
more useful for both theoretical and practical purposes.

We can treat larger values of $d$ as well, but at the cost of having a longer
list of possibilities.  For instance, with the hypotheses and notation
of the above result, if $d=11$ then $f$ permutes\/ $\F_q$ if and only if either $(*)$
holds or there is some $\psi\in\CC$ such that
\[
\left(\frac{\omega^{ike}-\omega^{-ike}}{\omega^{ie}-
\omega^{-ie}}\right)^{st} = \omega^{(2r+(k-1)vt)\psi(i)}
\]
for every $i\in\left(\F_{11}^*\right)^2$, where
$\CC$ is the union of the sets $\{mi: m\in\{\pm 3,\pm 5\}\}$,
$\{5m^3i^4 + m^7i^3 - 2mi^2 - 4m^5i: m\in\F_{11}^*\}$, and
$\{4m^3 i^4 + m^7i^3 - 2mi^2 - 5m^5i: m\in\F_{11}^*\}$.

\section{Preliminary lemma}

We begin with a simple lemma reducing the question whether a polynomial
permutes $\F_q$ to the question whether a related polynomial permutes
a particular subgroup of $\F_q^*$.  Here, for any positive integer $d$,
let $\mu_d$ denote the set of $d^{\operatorname{th}}$ roots of unity in
$\F_q$.

\begin{lemma}
\label{lwl}
Pick $d,r>0$ with $d\mid (q-1)$, and let $h\in\F_q[x]$.
Then $f(x):=x^r h(x^{(q-1)/d})$ permutes\/ $\F_q$ if and only if
both
\begin{enumerate}
\item $\gcd(r,(q-1)/d)=1$\, and
\item $x^r h(x)^{(q-1)/d}$ permutes $\mu_d$.
\end{enumerate}
\end{lemma}

\begin{proof}
Write $s:=(q-1)/d$.
For $\zeta\in\mu_s$, we have $f(\zeta x) = \zeta^r f(x)$.
Thus, if $f$ permutes $\F_q$ then $\gcd(r,s)=1$.
Conversely, if $\gcd(r,s)=1$ then the values of $f$ on $\F_q$
consist of all the $s^{\operatorname{th}}$ roots of the values
of
\[
f(x)^s = x^{rs} h(x^s)^s.
\]
But the values of $f(x)^s$ on $\F_q$ consist of $f(0)^s=0$
and the values of $g(x):=x^r h(x)^s$ on $(\F_q^*)^s$.
Thus, $f$ permutes $\F_q$ if and only if $g$ is bijective
on $(\F_q^*)^s=\mu_d$.
\end{proof}

\begin{remark}
A more complicated criterion for $f$ to permute $\F_q$ was given
by Wan and Lidl~\cite[Thm.\ 1.2]{WL}.
\end{remark}


\section{Proofs}

In this section we consider polynomials of the form $f(x):=x^r h_k(x^v)^t$,
where $h_k(x)=x^{k-1}+x^{k-2}+\dots+1$ and $r,v,k,t$ are positive integers.
We maintain this notation throughout this section, and we also define
$s:=\gcd(v,q-1)$, $d:=(q-1)/s$, and $e:=v/s$.  Note that $\gcd(d,e)=1$.

We begin with some easy cases with $d$ small:
\begin{prop}
\label{veryeasy}
If $d=1$ then $f(x)$ permutes $\F_q$ if and only if $\gcd(k,p)=\gcd(r,s)=1$.
If $d=2$ then $f(x)$ permutes $\F_q$ if and only if $\gcd(k,2p)=\gcd(r,s)=1$
and $k^{st}\equiv (-1)^{r+1}\pmod{p}$.
\end{prop}

\begin{proof}
By Lemma~\ref{lwl}, $f$ permutes $\F_q$ if and only if $\gcd(r,s)=1$ and
$g(x):=x^r h_k(x^e)^{st}$ permutes $\mu_d$.  If $d=1$, the latter condition
just says $\gcd(k,p)=1$, since $g(1)=k^{(q-1)t}$.  If $d=2$ then
we must have $h_k(-1)\ne 0$, so $k$ odd, whence $g(-1)=(-1)^r$; since
$g(1)=k^{st}$, the result follows.
\end{proof}

We could treat a few more values of $d$ by the same method as above, but
this requires handling several cases already for $d=3$.  We will return to
this question later in this section, after proving some results which simplify
the analysis.

Our next result gives necessary and sufficient conditions for $f$ to
permute $\F_q$; these conditions refine the ones we get directly from
Lemma~\ref{lwl}.

\begin{prop}
\label{pl}
$f$ permutes\/ $\F_q$ if and only if all of the following hold:
\begin{enumerate}
\item $\gcd(r,s)=\gcd(d,k)=1$
\item $\gcd(d,2r+vt(k-1))\le 2$
\item $k^{st}\equiv (-1)^{(d+1)(r+1)}\pmod{p}$
\item $g(x):=x^r((1-x^{ek})/(1-x^e))^{st}$ is injective on
$\mu_d\setminus\mu_1$
\item $(-1)^{(d+1)(r+1)} \notin g(\mu_d\setminus\mu_1)$.
\end{enumerate}
\end{prop}

\begin{proof}
By Lemma~\ref{lwl}, $f$ permutes $\F_q$ if and only if
$\gcd(r,s)=1$ and $\hat g(x):=x^r h_k(x^e)^{st}$ permutes $\mu_d$.
So assume $\gcd(r,s)=1$.  For $\zeta\in\mu_d\setminus\mu_1$, we have
\[
\hat g(\zeta)=\zeta^r \left(\frac{1 - \zeta^{ke}}{1-\zeta^e}\right)^{st},
\]
so $\hat g(\zeta)=0$ if and only if $\zeta\in\mu_{ke}$.
Thus, if $\hat g$ permutes $\mu_d$ then $\gcd(d,k)=1$ and
$\gcd(p,k)=1$ (since $\hat g(1)=k^{st}$).
Henceforth we assume $\gcd(pd,k)=1$, so $\hat g$ maps $\mu_d$ into $\mu_d$,
and thus bijectivity of $\hat g$ is equivalent to injectivity.

If $\hat g$ permutes $\mu_d$ then
\[
\prod_{\zeta\in\mu_d} \hat g(\zeta) = \prod_{\zeta\in\mu_d} \zeta = (-1)^{d+1},
\]
but we compute
\begin{align*}
\prod_{\zeta\in\mu_d} \hat g(\zeta) &= k^{st} 
\prod_{\zeta\in\mu_d\setminus\mu_1}
       \zeta^r \left(\frac{1 - \zeta^{ke}}{1-\zeta^e}\right)^{st} \\
&= (-1)^{(d+1)r} k^{st} \quad\text{ (since $\gcd(d,k)=1$).}
\end{align*}
Thus, if $\hat g$ permutes $\mu_d$ then $k^{st} = (-1)^{(d+1)(r+1)}$ in $\F_q$.

Next, if $\hat g$ permutes $\mu_d$ then for $\zeta\in\mu_d\setminus\mu_2$
we have $\hat g(\zeta)\ne \hat g(1/\zeta)$; but
$\hat g(1/\zeta)=\hat g(\zeta)/\zeta^{2r+set(k-1)}$,
so we conclude that $\gcd(d,2r+set(k-1))\le 2$.  The proof is complete.
\end{proof}

\begin{remark}
The fact that $k^{st}\equiv (-1)^{(d+1)(r+1)} \pmod{p}$ was
proved by Park and Lee \cite{PL} (in case $t=1$) by means of a lengthy
computation
of the determinants of some circulant matrices.  The case $t=1$ of
Proposition~\ref{pl}
improves the main result of \cite{AAW}; those authors gave some necessary
conditions for $f$ to permute $\F_q$, and some sufficient conditions, and
gave necessary and sufficient conditions in the special case that $d$ is
an odd prime less than $2p+1$.
\end{remark}

When $d$ is an odd prime, the criteria of Proposition~\ref{pl} can be
stated in terms of permutations of $\F_d$:

\begin{cor}
\label{neat}
Suppose the first three conditions of Proposition~\ref{pl} hold,
and $d$ is an odd prime.  Pick $\omega\in\F_q$ of order~$d$.
Then $f$ permutes $\F_q$ if and only if there exists $\theta\in\F_d[x]$
with $\theta(0)=0$ and $\deg(\theta)<(d-1)/2$ such that
$(2r+(k-1)vt)x+\theta(x^2)$ permutes $\F_d$
and, for every $i$ with $0<i<d/2$, we have
\[
\omega^{\theta(i^2)} = \left(\frac{\omega^{ike}-\omega^{-ike}}{\omega^{ie}-
\omega^{-ie}}\right)^{st}.
\]
\end{cor}

\begin{proof}
Since $d$ is odd, squaring permutes $\mu_d$, so condition (4) of
Proposition~\ref{pl} is equivalent to injectivity of $\hat g(x^2)$
on $\mu_d\setminus\mu_1$.  For $\zeta\in\mu_d\setminus\mu_1$, we have
\[
g(\zeta^2) = \zeta^{2r}\left(\frac{1-\zeta^{2ke}}{1-\zeta^{2e}}\right)^{st}
= \zeta^{2r+(k-1)est}\left(\frac{\zeta^{ke}-\zeta^{-ke}}{\zeta^e-\zeta^{-e}}
\right)^{st}.
\]
For $i\in\Z\setminus d\Z$, let $\psi(i)$ be the unique
element of $\Z/d\Z$ such that
\[
\omega^{\psi(i)} = \left(\frac{\omega^{ike}-\omega^{-ike}}{\omega^{ie}-
\omega^{-ie}}\right)^{st}.
\]
Defining $\psi(i)=0$ if $i\in d\Z$, it follows that $\psi$ induces
a map from $\Z/d\Z$ to itself, with the properties $\psi(-i)=\psi(i)$
and $g(\omega^{2i}) = \omega^{i(2r+(k-1)vt)+\psi(i)}$.  Conditions
(4) and (5)
are equivalent to bijectivity of the map $\chi:i\mapsto in+\psi(i)$ on
$\Z/d\Z$, where $n:=2r+(k-1)vt$.  Since $\psi(-i)=\psi(i)$, we can write
$\psi(i)=\theta(i^2)$
where $\theta\in\F_d[x]$ has degree less than $(d-1)/2$ and has no constant
term.
\end{proof}

For small $d$, there are only a few maps $\hat\theta:\F_d\to\F_d$ for which
$x+\hat\theta(x^2)$ permutes $\F_d$; this in turn yields manageable
descriptions of the possible permutation polynomials in these cases.  Assuming
$\hat\theta(0)=0$ and $\deg(\hat\theta)<(d-1)/2$, the only such map for $d=3$
and $d=5$ is $\hat\theta=0$.
For $d=7$ there are three possibilities for $\hat\theta$, namely
$\hat\theta=\mu x^2$ with $\mu\in\{0,2,-2\}$.  For $d=11$ there are $25$
possibilities for $\hat\theta$, but up to the equivalence $\hat\theta(x)\sim
\hat\theta(\alpha^2 x)/\alpha$ with $\alpha\in\F_d^*$, there are just five
possibilities.  For $d=13$ there are $133$ possibilities for $\hat\theta$,
including $14$ classes under the above equivalence.  We checked via computer
that, for these values of $d$, every such map $\hat\theta$ occurs as
$\theta/(2r+(k-1)vt)$ for some permutation polynomial $f$ as in
Corollary~\ref{neat}, even if we restrict to $k=2$ and $t=e=1$.

\begin{cor}
Suppose the first three conditions of Proposition~\ref{pl} hold,
and $d$ is an odd prime.  Pick $\omega\in\F_q$ of order~$d$.
\begin{enumerate}
\item[(a)] If
\begin{equation}
\tag{$*$}
\frac{\zeta^k-\zeta^{-k}}{\zeta-\zeta^{-1}}\in\mu_{st} \text{ for every
$\zeta\in\mu_d\setminus\mu_1$}
\end{equation}
then $f$ permutes\/ $\F_q$.
\item[(b)] If $d=3$ then $f$ always permutes\/ $\F_q$.
\item[(c)] If $d=5$ then $f$ permutes\/ $\F_q$ if and only if $(*)$ holds.
\item[(d)] If $d=7$ then $f$ permutes\/ $\F_q$ if and only if either $(*)$
holds or there exists $\epsilon\in\{1,-1\}$ such that
\[
\left(\frac{\omega^{ike}-\omega^{-ike}}{\omega^{ie}-
\omega^{-ie}}\right)^{st} = \omega^{2\epsilon (2r+(k-1)vt)i}
\]
for every $i\in\{1,2,4\}$.
\item[(e)] If $d=11$ then $f$ permutes\/ $\F_q$ if and only if either $(*)$
holds or there is some $\psi\in\CC$ such that
\[
\left(\frac{\omega^{ike}-\omega^{-ike}}{\omega^{ie}-
\omega^{-ie}}\right)^{st} = \omega^{(2r+(k-1)vt)\psi(i)}
\]
for every $i\in\left(\F_{11}^*\right)^2$, where
$\CC$ is the union of the sets $\{mi: m\in\{\pm 3,\pm 5\}\}$,
$\{5m^3i^4 + m^7i^3 - 2mi^2 - 4m^5i: m\in\F_{11}^*\}$, and
$\{4m^3 i^4 + m^7i^3 - 2mi^2 - 5m^5i: m\in\F_{11}^*\}$.
\end{enumerate}
\end{cor}

\begin{proof}
We maintain the notation of Corollary~\ref{neat}.
Condition $(*)$ is the trivial case $\theta=0$.
If $d=3$ or $d=5$, we plainly must have $\theta=0$ (as was first proved by
Betti in 1851 \cite{B}).
This proves the result for $d=5$.  For $d=3$,
condition (3) implies $k\equiv\pm 1\pmod{3}$, so for $\zeta\in\mu_d
\setminus\mu_1$ we have $\zeta^k-\zeta^{-k}=\pm(\zeta-
\zeta^{-1})$; since either $q$ or $s$ is even, this implies
$(\zeta^k-\zeta^{-k})^s=(\zeta-\zeta^{-1})^s$, so $(*)$ holds.

Suppose $d=7$, and write $n:=2r+(k-1)vt$; then $\gcd(7,n)=1$ by condition (2).
It is easy to determine the possibilities for $\theta$, as was first done
by Hermite in 1863 \cite{H}:
$\theta=\mu x^2$ where $\mu\in\{0,2n,-2n\}$.  The result follows.

The case $d=11$ is treated similarly.
\end{proof}



\end{document}